\documentclass{amsart}
\usepackage[english]{babel}
\usepackage{amsthm}
\usepackage{inputenc}
\usepackage{amssymb}
\usepackage{graphicx}
\newtheorem{thm}{Theorem}[section]
\newtheorem{lem}{Lemma}[section]
\newtheorem{pr}{Proposition}[section]

\numberwithin{equation}{section}

\begin{document}

\title{\textbf{On the solutions of the equation $x^3+ax=b$ in $\mathbb{Z}^*_3$
with coefficients from $\mathbb{Q}_3$}}

\author{Rikhsiboev I.M., Khudoyberdiyev A.Kh., Kurbanbaev T.K., Masutova K.K.}

\date{}
\maketitle

\begin{abstract}

Recall that in \cite{Mukh} it is obtained the criteria solvability
of the equation $x^3+ax=b$ in $\mathbb{Z}^*_p,$ $\mathbb{Z}_p$ and
$\mathbb{Q}_p$ for $p>3$. Since any $p$-adic number $x$ has a
unique form $x=p^kx^*,$ where $x^* \in \mathbb{Z}^*_p $ and $k \in
\mathbb{Z},$ in \cite{Mukh} it is also shown that from the
criteria in $\mathbb{Z}^*_p$ it follows the criteria in
$\mathbb{Z}_p$ and $\mathbb{Q}_p$. In this paper we provide the
algorithm of finding the solutions of the equation $x^3+ax=b$ in
$\mathbb{Z}^*_3$ with coefficients from $\mathbb{Q}_3$.
\end{abstract}



\section{Introduction}

In the present time description of different structures in
mathematics are studying over field of $p$-adic numbers. In
particular, $p$-adic analysis is one of the intensive
developing directions of modern mathematics. Numerous applications of
$p$-adic numbers have found their own reflection in the theory of
$p$-adic differential equations, $p$-adic theory of probabilities,
$p$-adic mathematical physics, algebras over $p$-adic numbers and
others.

The field of $p$-adic numbers were introduced by German mathematician
K. Hensel at the end of the 19th century \cite{Hen}. The investigation of $p$-adic numbers were motivated primarily by an attempt
to bring the ideas and techniques of the power series into number theory.
Their canonical representation is similar to expansion of analytical functions
in power series, which is analogy between algebraic numbers and algebraic functions.
There are several books devoted to study
$p$-adic numbers and $p$-adic analysis \cite{Buh}, \cite{Katok}, \cite{Kobl}, \cite{Vlad}.

Classification of algebras in small dimensions plays important
role for the studying of properties of varieties of algebras. It
is known that the problem of classification of finite dimensional
algebras involves a study on equations for structural
constants, i.e. to the decision of some systems of the equations
in the corresponding field. Classifications of complex Leibniz
algebras have been investigated in \cite{Alb1-lengh} \cite{AO1}, \cite{Omir3},
\cite{RS1} and many other works. In similar complex case, the problem of
classification in $p$-adic case is reduced to the solution of the
equations in the field. The classifications of Leibniz algebras over the field of $p$-adic numbers
have been obtained in \cite{Ayup},\cite{Khud}, \cite{Lad}.

In the field of complex numbers the fundamental
Abel's theorem about insolvability in radicals of general equation
of $n$-th degree $(n>5)$ is well known. In this field square equation is solved
by discriminant, for cubic equation Cardano's
formulas were widely applied. In the field of $p$-adic numbers square equation does not
always has a solution. Note that the criteria of solvability of
the equation $x^2=a$ is given in \cite{Vlad} and in \cite{Cas},
\cite{Mukh2} we can find the solvability criteria for the equation
$x^q =a,$ where $q$ is an arbitrary natural number.

In this paper we consider $p-$adic cubic equation
$y^3+ry^2+sy+t=0.$  By replacing $y=x-\displaystyle\frac{r}{3},$
this equation become the so-called depressed cubic equation
\begin{equation}\label{eq1.1}x^3+ax=b.\end{equation}

The solvability criterion for the cubic equation $x^3+ax=b $ over $3$-adic numbers is different from the case $p>3.$
Note that solvability criteria for $p>3$ is obtained in \cite{Mukh} and for
$p=3$ with condition of $|a|_3\neq \frac 1 3$ is given in \cite{Kur}.

In this paper by using the results of \cite{Kur}, the algorithm of finding the solutions of
the equation $x^3+ax=b$ in $\mathbb{Z}^*_3$ with coefficients from
$\mathbb{Q}_3$ will be provided.

\section{Preliminaries}

Let $\mathbb{Q}$ be a field of rational numbers. Every rational
number $x\neq0$ can be represented by the form
$x=p^{\gamma(x)}\displaystyle\frac{n}{m}$, where $n,\gamma(x)\in
\mathbb{Z},$ $m$ is a positive integer, $(p,n)=1,\ (p,m)=1$ and
$p$ is a fixed prime number. In $\mathbb{Q}$ a norm has been defined as follows:
$$|x|_p=\left\{\begin{array}{lll}
p^{-\gamma(x)}, &  x\neq 0, \\
0, & x=0. \\
\end{array}\right.$$

The norm $|x|_p$ is called a \emph{$p$-adic norm} of $x$ and it
satisfies so called the strong triangle inequality. The completion
of $\mathbb{Q}$ with respect to $p$-adic norm defines the $p$-adic
field which is denoted by $\mathbb{Q}_p$ (\cite{Katok},\cite{Vlad}). It is well
known that any $p$-adic number $x\neq0$ can be uniquely
represented in the canonical form
$$x=p^{\gamma(x)}(x_0+x_1p+x_2p^2+\dots),$$ where
$\gamma=\gamma(x)\in\mathbb{Z}$  and  $x_j$ are integers, $0\leq
x_j \leq p-1$, $x_0\neq 0$, ($j=0,1,\dots$). $p$-Adic number $x,$
for which $|x|_p\leq1$, is called \emph{integer} $p$-adic number,
and the set of such numbers is denoted by $\mathbb{Z}_p.$ Integer
$x\in\mathbb{Z}_p,$ for which $|x|_p=1,$ is called \emph{unit} of
$\mathbb{Z}_p,$ and their set is denoted by $\mathbb{Z}^*_p$.

For any numbers $a$ and $m$ it is known the following result.

\begin{thm}\label{t1} \cite{Buh}. If $(a,m)=1,$ then a
congruence $ax\equiv b \, (mod \, m)$ has one and only one
solution.
\end{thm}

We also need the following Lemma.

\begin{lem}\label{l1} \cite{Cas}. The following is true:
$$\left(\sum_{i=0}^\infty
x_ip^i\right)^q=x_0^q+\sum_{k=1}^\infty\left(qx_0^{q-1}x_k+N_k(x_0,x_1,\dots,x_{k-1})\right)p^k,$$
where $x_0\ne 0$, $0\leq x_j\leq p-1$, $N_1=0$ and for $k\geq 2$

$$N_k=N_k(x_0,\dots,x_{k-1})=
\sum\limits_{{m_0,m_1,\dots,m_{k-1}:\atop \sum_{i=0}^{k-1}m_i=q,\
\ \sum_{i=1}^{k-1}im_i=k}} {q!\over m_0!m_1!\dots
m_{k-1}!}x_0^{m_0}x_1^{m_1}\dots x_{k-1}^{m_{k-1}}.$$
\end{lem}

From the lemma \ref{l1} by $q=3$ we have

$$\left(\sum_{i=0}^\infty
x_ip^i\right)^3=x_0^3+\sum_{k=1}^\infty\left(3x_0^2x_k+N_k(x_0,x_1,\dots,x_{k-1})\right)p^k.$$

For $j\leq k$ we put
$$P_k^j=P_k^j(x_0,x_1,\dots, x_{j-1})=\sum\limits_{{m_0,m_1,\dots,m_{j-1}:\atop
\sum_{i=0}^{j-1}m_i=3,\ \ \sum_{i=1}^{j-1}im_i=k}} {6\over
m_0!m_1!\dots m_{j-1}!}x_0^{m_0}x_1^{m_1}\dots
x_{j-1}^{m_{j-1}}.$$

Also the following identity is true:
\begin{equation}\label{eq2.1}\left(\sum_{i=0}^\infty
a_ip^i\right)\left(\sum_{j=0}^\infty x_jp^j\right)=
\sum_{k=0}^\infty \left(\sum_{s=0}^k x_s a_{k-s}\right)p^k.
\end{equation}

\

\

\section{The main result}

In this paper we study the cubic equation \eqref{eq1.1} over the
field $3$-adic numbers, i.e. $a,b, x \in \mathbb{Q}_3.$

Put
$$x=3^{\gamma(x)}(x_0+x_13+\dots),\quad a=3^{\gamma(a)}(a_0+a_13+\dots),\quad
b=3^{\gamma(b)}(b_0+b_13+\dots),$$ where $x_j, a_j, b_j
\in\{0,1,2\}, x_0, a_0, b_0\neq0, (j=0,1,\dots)$.

Since any $3$-adic number $x$ has a unique form $x=3^kx^*,$ where $x^* \in \mathbb{Z}^*_3$
and $k \in \mathbb{Z},$ we will be limited to search a decision from $\mathbb{Z}^*_3,$ i.e. $\gamma(x)=0$.

Putting the canonical form of $a,b$ and $x$ in \eqref{eq1.1}, we
get
$$\left(\sum_{k=0}^\infty
x_k 3^k\right)^3 + 3^{\gamma(a)}\left(\sum_{k=0}^\infty a_k
3^k\right)\left(\sum_{k=0}^\infty x_k 3^k\right)= 3^{\gamma(b)}
\sum_{k=0}^\infty b_k 3^k.$$

By Lemma \ref{l1} and equality \eqref{eq2.1}, the equation
\eqref{eq1.1} becomes to the following form:
\begin{equation}\label{eq3.1}
\begin{array}{l}x_0^3+\sum\limits_{k=1}^\infty\left(3x_0^2x_k+N_k(x_0,x_1,\dots,x_{k-1})\right)3^k
+3^{\gamma(a)} \left(a_0x_0+\sum\limits_{k=1}^\infty
\left(\sum\limits_{s=0}^k x_s a_{k-s}\right)3^k\right)\\
=3^{\gamma(b)}\left(b_0+ \sum\limits_{k=1}^\infty b_k
3^k\right).\end{array}\end{equation}

\begin{pr} \label{p1} If one of the following conditions:
$$\begin{array}{llllllll} 1) & \gamma(a)=0 & \mbox{\emph{and}} &
\gamma(b)<0;&
2) & \gamma(a)>0 & \mbox{\emph{and}} & \gamma(b)>0;\\[1mm]
3) & \gamma(a)>0 & \mbox{\emph{and}} & \gamma(b)<0;&
4) & \gamma(a)<0 & \mbox{\emph{and}} & \gamma(b)=0;\\[1mm]
5) & \gamma(a)<0 & \mbox{\emph{and}} & \gamma(b)>0,\end{array}$$
is fulfilled, then the equation \eqref{eq1.1} has not a solution
in $\mathbb{Z}^*_3.$
\end{pr}

\begin{proof} 1) Let $\gamma(a)=0$ and $\gamma(b)<0.$
Multiplying the equation \eqref{eq3.1} by $3^{-\gamma(b)},$ we get
the following congruence $b_0 \equiv 0 \,\ (mod \, 3),$ which is
not correct. Consequently, the equation \eqref{eq1.1} has no
solution in $\mathbb{Z}^*_3.$

2) Let $\gamma(a)>0$ and $\gamma(b)>0.$ Then from \eqref{eq3.1} it
follows a congruence $x^3_0\equiv 0 \,\ (mod \, 3),$ which has no
a nonzero solution.Therefore, in $\mathbb{Z}^*_3$ the equation
\eqref{eq1.1} does not have a solution.

In other cases we analogously get the congruences
$$b_0 \equiv 0 \,\ (mod \, 3)\quad\mbox{or} \quad a_0x_0\equiv 0 \,\ (mod \, 3),$$
which are not hold. Therefore, in $\mathbb{Z}^*_3$ there is no solution.
\end{proof}

From the Proposition \ref{p1} we have that the cubic equation may have a
solution if one of the following four cases
$$\begin{array}{cccccc} 1) & \gamma(a)=0,& \gamma(b)=0,  &  2) &  \gamma(a)=0,& \gamma(b)>0,\\[1mm]
3) &  \gamma(a)<0, & \gamma(b)<0, &   4) &  \gamma(a)>0,
&\gamma(b)=0 \end{array}$$ is hold.

In the following theorem we present an algorithm of finding of the
equation $x^3+ax=b$ for the first case.

\begin{thm}\label{t2} Let $\gamma(a)=\gamma(b)=0$ and
$a_0=1.$ Then $x$ to be a solution of the equation \eqref{eq1.1}
in $\mathbb{Z}^*_3$ if and only if the congruences
$$x^3_0+a_0x_0\equiv b_0 \,\ (mod \, 3),$$
$$x_1a_0+x_0a_1+N_1(x_0)+M_1(x_0)\equiv b_1 \ (mod \, 3),$$
$$x_ka_0+ \dots +
x_0a_k+x_0^2x_{k-1}+N_k(x_0,x_1,\dots,x_{k-1})+M_k(x_0,x_1,\dots,x_{k-1})\equiv
b_k \ (mod \, 3), \ k \geq 2$$ are fulfilled, where integers
$M_{k}(x_0,\dots,x_{k-1})$ are defined consequently from the
following correlations
$$x^3_0+a_0x_0 = b_0 +  M_1(x_0)\cdot 3,$$
$$x_1a_0+x_0a_1+N_1(x_0)=  b_1 - M_1(x_0) + M_2(x_0,x_1)\cdot3,$$
$$x_{k-1}a_0+ x_{k-2}a_1+ \dots + x_0a_{k-1}+x_0^2x_{k-2}+
N_{k-1}(x_0,x_1,\dots,x_{k-2})=$$$$= b_{k-1}
-M_{k-1}(x_0,x_1,\dots,x_{k-2})+ M_k(x_0,x_1,\dots,x_{k-1}) \cdot
3, \quad k \geq 3.$$
\end{thm}

\begin{proof}
Let $$x_0+x_1\cdot 3+x_2\cdot 3^2+\dots,\,\ 0\leq x_j \leq 2, \,\
x_0\neq 0, (j=0,1,\dots)$$ is a solution of the equation
\eqref{eq1.1}, then equality \eqref{eq3.1} becomes
$$x_0^3+\sum\limits_{k=1}^\infty\left(3x_0^2x_k+N_k(x_0,x_1,\dots,x_{k-1})\right)3^k+
a_0x_0+\sum\limits_{k=1}^\infty
\left(\sum\limits_{s=0}^k x_s
a_{k-s}\right)3^k=b_0+
\sum\limits_{k=1}^\infty b_k 3^k.$$

So we have
$$x^3_0+a_0x_0 + (x_1a_0+x_0a_1+N_1(x_0))\cdot 3 +$$
$$+ \sum\limits_{k=2}^\infty
\left(x_ka_0+ x_{k-1}a_1+ \dots +
x_0a_k+x_0^2x_{k-1}+N_k(x_0,x_1,\dots,x_{k-1})\right)3^k = b_0+
\sum\limits_{k=1}^\infty b_k 3^k,$$ from which it follows
the necessity in fulfilling the congruences of the theorem.


Now let $x$ is satisfied the congruences of the theorem. Since
$(a_0,3) =1,$ then by Theorem \ref{t1} it implies that these
congruences have the solutions $x_k.$

Then
$$x_0^3+\sum\limits_{k=1}^\infty\left(3x_0^2x_k+N_k(x_0,x_1,\dots,x_{k-1})\right)3^k
+ a_0x_0+\sum\limits_{k=1}^\infty \left(\sum\limits_{s=0}^k x_s
a_{k-s}\right)3^k=$$
$$=x_0^3+a_0x_0 + (N_1+x_0a_1+a_0x_1)3+$$
$$+\sum\limits_{k=2}^\infty\left(x_0^2x_{k-1}+N_k(x_0,x_1,\dots,x_{k-1})+x_0a_k+x_1a_{k-1}+ \dots +x_{k-1}a_1+x_ka_0\right)3^k=$$
$$=b_0 +M_1(x_0)\cdot3+(b_1 - M_1(x_0)+ M_2(x_0,x_1)\cdot3)\cdot3+$$
$$+\sum\limits_{k=2}^\infty\left(b_k - M_k(x_0,x_1,\dots,x_{k-1})+
M_{k+1}(x_0,x_1,\dots,x_k)\cdot3\right)\cdot3^k= b_0+
\sum\limits_{k=1}^\infty b_k 3^k.$$

Therefore, we show that $x=\sum\limits_{k=0}^\infty x_k 3^k$ is  a
solution of the equation \eqref{eq1.1}. \end{proof}

Let us examine a case $\gamma(a)=0, \gamma(b)>0$ and get necessary
and sufficient conditions for a solution of the equation
\eqref{eq1.1}.

\begin{thm}\label{t3} Let $\gamma(a)=0,$ $\gamma(b)=m>0$ and
$a_0=2.$ Then $x$ to be a solution of the equation \eqref{eq1.1}
in $\mathbb{Z}^*_3$ if and only if the congruences
$$x^3_0+a_0x_0\equiv 0 \, (mod \, 3),$$
$$x_1a_0+x_0a_1+N_1(x_0)+M_1(x_0)\equiv 0 \ (mod \, 3),$$
$$x_ka_0+  \dots +
x_0a_k+x_0^2x_{k-1}+N_k(x_0,\dots,x_{k-1})+M_k(x_0,\dots,x_{k-1})\equiv
0 \ (mod \, 3),\,  2\leq k \leq m-1,$$
$$x_ka_0+  \dots +
x_0a_k+x_0^2x_{k-1}+N_k(x_0,\dots,x_{k-1})+M_k(x_0,\dots,x_{k-1})\equiv
b_{k-m} \ (mod \, 3), \, k \geq m$$ are fulfilled, where integers
$M_{k}(x_0,x_1,\dots,x_{k-1})$ are defined consequently from the
following correlations
$$x^3_0+2x_0 = M_1(x_0)\cdot3,$$
$$x_1a_0+x_0a_1+N_1(x_0)= -M_1(x_0)+3M_2(x_0,x_1),$$
$$x_ka_0+ x_{k-1}a_1+ \dots +
x_0a_k+x_0^2x_{k-1}+N_k(x_0,x_1,\dots,x_{k-1})=$$
$$=-M_k(x_0,x_1,\dots,x_{k-1}) +3M_{k+1}(x_0,x_1,\dots,x_k), \, 2\leq k \leq
m-1,$$
$$x_ka_0+ x_{k-1}a_1+ \dots +
x_0a_k+x_0^2x_{k-1}+N_k(x_0,x_1,\dots,x_{k-1})=$$$$ =b_{k-m}-
M_k(x_0,x_1,\dots,x_{k-1}) +  3M_{k+1}(x_0,x_1,\dots,x_k), \, k
\geq m.$$ \end{thm}

\begin{proof} Let $x$ is a solution of the equation
\eqref{eq1.1}, then equality \eqref{eq3.1} becomes
$$x_0^3+\sum\limits_{k=1}^\infty\left(3x_0^2x_k+N_k(x_0,x_1,\dots,x_{k-1})\right)3^k+
 a_0x_0+\sum\limits_{k=1}^\infty
\left(\sum\limits_{s=0}^k x_s a_{k-s}\right)3^k=3^m\left(b_0+
\sum\limits_{k=1}^\infty b_k 3^k\right).$$

Therefore, we have
$$x^3_0+a_0x_0 + (x_1a_0+x_0a_1+N_1(x_0))\cdot 3 +$$
$$+ \sum\limits_{k=2}^\infty
\left(x_ka_0+ x_{k-1}a_1+ \dots +
x_0a_k+x_0^2x_{k-1}+N_k(x_0,\dots,x_{k-1})\right)3^k
=3^m\left(b_0+ \sum\limits_{k=1}^\infty b_k 3^k\right),$$ from
which it follows the necessity in fulfilling the congruences of
the theorem.

Now let $x$ is satisfied the congruences of the theorem. Since
$(a_0,3) =1,$ then by Theorem \ref{t1} there are solutions $x_k$
of the congruences.

Putting element $x$ to the equality \eqref{eq3.1}, we have
$$x_0^3+\sum\limits_{k=1}^\infty\left(3x_0^2x_k+N_k(x_0,x_1,\dots,x_{k-1})\right)3^k
+ a_0x_0+\sum\limits_{k=1}^\infty \left(\sum\limits_{s=0}^k x_s
a_{k-s}\right)3^k=$$
$$=x_0^3+a_0x_0 + (N_1+x_0a_1+a_0x_1)3+$$
$$+\sum\limits_{k=2}^\infty\left(x_0^2x_{k-1}+N_k(x_0,x_1,\dots,x_{k-1})
+x_0a_k+x_1a_{k-1}+ \dots +x_{k-1}a_1+x_ka_0\right)3^k=$$
$$=M_1(x_0)\cdot3+(- M_1(x_0)+ M_2(x_0,x_1)\cdot3)\cdot3+$$
$$+\sum\limits_{k=2}^{m-1}\left(- M_k(x_0,x_1,\dots,x_{k-1})+
M_{k+1}(x_0,x_1,\dots,x_k)\cdot3\right)\cdot3^k+$$
$$+\sum\limits_{k=m}^\infty\left(b_{k-m} - M_k(x_0,x_1,\dots,x_{k-1})+
M_{k+1}(x_0,x_1,\dots,x_k)\cdot3\right)\cdot3^k= 3^m\left(b_0+
\sum\limits_{k=1}^\infty b_k 3^k\right).$$

Therefore, we show that $x$ is  a solution of the equation
\eqref{eq1.1}.
\end{proof}

The following theorem gives necessary and sufficient conditions
for a solution of the equation \eqref{eq1.1} for the case
$\gamma(a)<0$ and $\gamma(b)<0.$

\begin{thm}\label{t4} Let $$\gamma(a)=\gamma(b)=-m<0,\quad (m>0).$$ Then $x$ to be a solution of the equation
\eqref{eq1.1} in  $\mathbb{Z}^*_3$ if and only if the next
congruences
$$a_0x_0\equiv b_0 \ (mod \, 3),$$
$$x_ka_0+ x_{k-1}a_1+ \dots + x_0a_k
+M_k(x_0,x_1,\dots,x_{k-1})\equiv b_k \ (mod \, 3),\, 1\leq k \leq
m-1,$$
$$x_ma_0+ x_{m-1}a_1+ \dots + x_0a_m+x_0^3
+M_m(x_0,x_1,\dots,x_{m-1})\equiv b_m \ (mod \, 3),$$
$$x_{m+1}a_0+ x_ma_1+\dots+x_0a_{m+1}
+M_{m+1}(x_0,x_1,\dots,x_m)\equiv b_{m+1} \ (mod \, 3),$$
$$x_ka_0+ x_{k-1}a_1+ \dots +
x_0a_k+x_0^2x_{k-m-1}+N_{k-m}(x_0,x_1,\dots,x_{k-m-1})+$$
$$+M_k(x_0,x_1,\dots,x_{k-1})\equiv b_k \ (mod \, 3),\,k \geq
m+2$$ are fulfilled, where integers $M_{k}(x_0,x_1,\dots,x_{k-1})$
are defined consequently from the equalities
$$a_0x_0 =  b_0 +M_1(x_0)\cdot3,$$
$$x_ka_0+x_{k-1}a_1+ \dots + x_0a_k = b_k-M_k(x_0,\dots,x_{k-1})
+3M_{k+1}(x_0,\dots,x_k),\,1\leq k \leq m-1,$$
$$x_ma_0+ x_{m-1}a_1+ \dots +x_0a_m
+x_0^3=b_m-M_m(x_0,x_1,\dots,x_{m-1})+3M_{m+1}(x_0,x_1,\dots,x_m),$$
$$x_{m+1}a_0+ x_ma_1+\dots+x_0a_{m+1}=b_{m+1}-
M_{m+1}(x_0,\dots,x_m)+3M_{m+2}(x_0,\dots,x_{m+1}),$$
$$x_ka_0+ x_{k-1}a_1+ \dots +x_0a_k+x_0^2x_{k-m-1}+
N_{k-m}(x_0,x_1,\dots,x_{k-m-1})=$$
$$=b_k-M_k(x_0,x_1,\dots,x_{k-1})+3M_{k+1}(x_0,x_1,\dots,x_k),\,k\geq
m+2.$$ \end{thm}

\begin{proof} The proof of the Theorem can be obtained by similar way to the proofs of
Theorems \ref{t2}, \ref{t3}.

\end{proof}

Examining various cases of $\gamma(a)$ and $\gamma(b)$ we need to
study only the case $\gamma(a)>0$ and $\gamma(b)=0.$ Because of
appearance of uncertainty of a solution, we divide this case to
$\gamma(a)>1$ and $\gamma(a)=1.$

\begin{thm}\label{t5} Let $\gamma(a)= 2,$ $\gamma(b)=0$
and $(b_0, b_1) = (1, 0)$  or $(2, 2).$ Then $x$ to be a solution
of the equation \eqref{eq1.1} \emph{in} $\mathbb{Z}^*_3$ if and
only if he next congruences
$$x^3_0\equiv b_0(mod \, 3),$$
$$x^3_0\equiv b_0 +b_1\cdot 3\, (mod \, 9),$$
$$x_0^2x_1+x_0a_0 +M_1(x_0)\equiv b_2(mod \ 3),$$
$$x_0^2x_2+P_3^2(x_0,x_1)+x_1a_0+x_0a_1+x_0x_1^2+M_2(x_0,x_1)\equiv b_3(mod \ 3),$$
$$x_0^2x_{k-1}+P_k^{k-1}(x_0,x_1,\dots,x_{k-2})+2x_0x_1x_{k-2}+x_{k-2}a_0
+ x_{k-3}a_1+ \dots + x_0a_{k-2}+$$$$+M_{k-1}(x_0,x_1,\dots,
x_{k-2})\equiv b_k(mod \ 3), \quad k\geq 4$$ are fulfilled, where
integers $M_{k}(x_0,x_1,\dots,x_{k-1})$ are defined from the
equalities
$$x^3_0 = b_0 +b_1\cdot 3 + M_1(x_0)\cdot 9,$$
$$x_0^2x_1+x_0a_0 = b_2 - M_1(x_0)+3M_2(x_0,x_1),$$
$$x_0^2x_2+P_3^2(x_0,x_1)+x_1a_0+x_0a_1+x_0x_1^2 = b_3-M_2(x_0,x_1)+3M_3(x_0,x_1,x_2),
$$$$x_0^2x_{k-1}+P_k^{k-1}(x_0,x_1,\dots,x_{k-2})+2x_0x_1x_{k-2}+x_{k-2}a_0
+ x_{k-3}a_1+ \dots + x_0a_{k-2}=$$$$=b_k-M_{k-1}(x_0,x_1,\dots,
x_{k-2})+3M_k(x_0,x_1,\dots, x_{k-1}), \quad k\geq 4.$$ \end{thm}

\begin{proof} Analogously to the proof of the Theorem \ref{t2}.

\end{proof}

\begin{thm}\label{t6} Let $\gamma(a)= 3,$ $\gamma(b)=0$ and $(b_0,
b_1) = (1, 0)$ or $(2, 2).$ Then $x$ to be a solution of the
equation \eqref{eq1.1} in  $\mathbb{Z}^*_3$ if and only if he next
congruences
$$x^3_0\equiv b_0(mod \, 3),$$
$$x^3_0\equiv b_0 +b_1\cdot 3\, (mod \, 9),$$
$$x_0^2x_1 +M_1(x_0)\equiv b_2(mod \ 3),$$
$$x_0^2x_2+P_3^2(x_0,x_1)+x_0a_0+x_0x_1^2+M_2(x_0,x_1)\equiv b_3(mod
\ 3),$$
$$x_0^2x_{k-1}+P_k^{k-1}(x_0,x_1,\dots,x_{k-2})+2x_0x_1x_{k-2}+x_{k-3}a_0
+ x_{k-4}a_1+ \dots + x_0a_{k-3}+$$$$+M_{k-1}(x_0,x_1,\dots,
x_{k-2})\equiv b_k(mod \ 3), \quad k\geq 4$$ are fulfilled, where
integers $M_{k}(x_0,x_1,\dots,x_{k-1})$ are defined from the
equalities
$$x^3_0 = b_0 +b_1\cdot 3 + M_1(x_0)\cdot 9,$$
$$x_0^2x_1 = b_2 - M_1(x_0)+3M_2(x_0,x_1),$$
$$x_0^2x_2+P_3^2(x_0,x_1)+x_0a_0+x_0x_1^2 = b_3-M_2(x_0,x_1)+3M_3(x_0,x_1,x_2),
$$$$x_0^2x_{k-1}+P_k^{k-1}(x_0,x_1,\dots,x_{k-2})+2x_0x_1x_{k-2}+x_{k-3}a_0
+ x_{k-4}a_1+ \dots + x_0a_{k-3}=$$$$=b_k-M_{k-1}(x_0,x_1,\dots,
x_{k-2})+3M_k(x_0,x_1,\dots, x_{k-1}), \quad k\geq 4.$$ \end{thm}

\begin{proof} Analogously to the proof of the Theorem \ref{t2}. \end{proof}

Similarly to the Theorem \ref{t5}, it is proved the following

\begin{thm}\label{t7} Let $\gamma(a)= m \geq 4,$
$\gamma(b)=0$ and $(b_0, b_1) = (1, 0)$  or $(2, 2).$ Then $x$ to
be a solution of the equation \eqref{eq1.1} \emph{in}
$\mathbb{Z}^*_3$ if and only if he next congruences
$$x^3_0\equiv b_0(mod \, 3),$$
$$x^3_0\equiv b_0 +b_1\cdot 3\, (mod \, 9),$$
$$x_0^2x_1 +M_1(x_0)\equiv b_2(mod \ 3),$$
$$x_0^2x_2+P_3^2(x_0,x_1)+x_0x_1^2+M_2(x_0,x_1)\equiv b_3(mod
\ 3),$$
$$x_0^2x_{k-1}+P_k^{k-1}(x_0,x_1,\dots,x_{k-2})+2x_0x_1x_{k-2}\equiv b_k(mod \ 3), \quad 4\leq k \leq {m-1},$$
$$x_0^2x_{m-1}+P_m^{m-1}(x_0,x_1,\dots,x_{m-2})+2x_0x_1x_{m-2}+x_0a_0\equiv b_m(mod \ 3),$$
$$x_0^2x_{k-1}+P_k^{k-1}(x_0,x_1,\dots,x_{k-2})+2x_0x_1x_{k-2}+x_{k-m}a_0
+ \dots + x_0a_{k-m}\equiv b_k(mod \ 3), \quad k\geq{m+1}$$ are
fulfilled, where integers $M_{k}(x_0,x_1,\dots,x_{k-1})$ are
defined from the equalities
$$x^3_0 = b_0 +b_1\cdot 3 + M_1(x_0)\cdot 9,$$
$$x_0^2x_1 = b_2 - M_1(x_0)+3M_2(x_0,x_1),$$
$$x_0^2x_2+P_3^2(x_0,x_1)+x_0x_1^2 = b_3-M_2(x_0,x_1)+3M_3(x_0,x_1,x_2),
$$$$x_0^2x_{k-1}+P_k^{k-1}(x_0,x_1,\dots,x_{k-2})+2x_0x_1x_{k-2}=$$
$$=b_k-M_{k-1}(x_0,x_1,\dots,
x_{k-2})+3M_k(x_0,x_1,\dots, x_{k-1}), \quad 4\leq k \leq {m-1},$$
$$x_0^2x_{m-1}+P_m^{m-1}(x_0,\dots,x_{m-2})+2x_0x_1x_{m-2}+x_0a_0=$$$$= b_m-M_{m-1}(x_0,\dots,
x_{m-2})+3M_m(x_0,\dots, x_{m-1}),$$
$$x_0^2x_{k-1}+P_k^{k-1}(x_0,x_1,\dots,x_{k-2})+2x_0x_1x_{k-2}+x_{k-m}a_0
+ \dots + x_0a_{k-m}=$$$$=b_k-M_{k-1}(x_0,x_1,\dots,
x_{k-2})+3M_k(x_0,x_1,\dots, x_{k-1}), \quad k\geq m+1.$$
\end{thm}

Now we consider the equality \eqref{eq3.1} with $\gamma(a)=1,$
$\gamma(b)=0.$ Put
$$A_0=x^2_0+a_0,\quad  A_k=\displaystyle\frac{A_{k-1}}{3}+a_k+R_k, \quad
\mbox{where} \quad  R_k=\sum\limits_{j=0}^{k}x_j x_{k-j}, \,
k\geq1,$$
$$N'_j=\left\{\begin{array}{ll}
\displaystyle\frac{N_{j-1}}{3}, \, \mbox{если} \, j=3s-1, \\
[1mm]
\displaystyle\frac{N_{j-1}}{3}+x^3_{\frac j 3}, \, \mbox{если} \, j=3s, \,  \\
[1mm]
\displaystyle\frac{N_{j-1}-x^3_{\frac {j-1} 3}}{3}, \, \mbox{если} \, j=3s+1, \\
[1mm]
\end{array}\right. \quad S_j^i=\left\{\begin{array}{ll}
\displaystyle\frac{P_{j-1}^i}{3}, \, \mbox{если} \, j=3s-1, \\
[1mm]
\displaystyle\frac{P_{j-1}^i}{3}+x^3_{\frac j 3}, \, \mbox{если} \, j=3s, \,  \\
[1mm]
\displaystyle\frac{P_{j-1}^i-x^3_{\frac {j-1} 3}}{3}, \, \mbox{если} \, j=3s+1. \\
[1mm]
\end{array}\right.$$

\begin{thm}\label{t8} Let $\gamma(a)=1,$ $\gamma(b)=0$ and $x\in
\mathbb{Z}^{*}_3$ to be so that $A_0=x^2_0+a_0\not\equiv 0 \, (mod
\,3).$ Then $x$ to be a solution of the equation \eqref{eq1.1} in
$\mathbb{Z}^*_3$ if and only if the congruences
$$x^3_0\equiv b_0 \, (mod \, 3),$$
$$x_0a_0 + M_1(x_0)\equiv b_1 \, (mod \, 3),$$
$$(x_0^2+a_0)x_{k-1}+x_{k-2}a_1 + \dots +
x_0a_{k-1}+N'_k(x_0,\dots,x_{k-2})+M_k(x_0,\dots,x_{k-2}) \equiv
b_k\ (mod \ 3), \ k \geq 2$$ are faithfully, where
$M_1(x_0)=\frac{x^3_0-b_0} {3}$ and integers
$M_{k}(x_0,\dots,x_{k-2}), (k \geq 2)$ are defined from the
equalities
$$ x_0a_0 + M_1(x_0) = b_1 + M_2(x_0)\cdot 3,$$
$$(x_0^2+a_0)x_{k-1}+N'_k(x_0,x_1,\dots,x_{k-2})+x_{k-2}a_1 + \dots +
x_0a_{k-1}+$$
$$+M_k(x_0,\dots,x_{k-2}) =  b_k + M_{k+1}(x_0,\dots,x_{k-1})\cdot 3, \quad k \geq 2.$$
\end{thm}

\begin{proof}

Let the congruences $x^3_0\equiv b_0 \, (mod \, 3),\ x_0a_0 +
\frac{x^3_0-b_0} {3}\equiv b_1 \, (mod \, 3),$ has a solution
$x_0.$ Then denote by $M_2(x_0)$ the number satisfying the
equality $3M_2(x_0) = x_0a_0 + M_1(x_0) - b_1.$

Using Theorem \ref{t1}, we have existence of solutions $x_k$ of
the congruences
$$(x_0^2+a_0)x_{k-1}+x_{k-2}a_1 + \dots +
x_0a_{k-1}+N'_k(x_0,\dots,x_{k-2})+M_k(x_0,\dots,x_{k-2}) \equiv
b_k\ (mod \ 3), \quad k \geq 2.$$

The next chain of equalities
$$x_0^3+a_0x_0\cdot 3 + \sum\limits_{k=2}^\infty\left((x_0^2+a_0)x_{k-1}+N'_k(x_0,x_1,\dots,x_{k-2})+x_{k-2}a_1 + \dots + x_0a_{k-1}\right)3^k=$$
$$=b_0 +M_1(x_0)\cdot 3+(b_1 - M_1(x_0)+ M_2(x_0)\cdot 3)\cdot 3+$$
$$+\sum\limits_{k=2}^\infty\left(b_k - M_k(x_0,x_1,\dots,x_{k-2})+
M_{k+1}(x_0,x_1,\dots,x_{k-1})\cdot 3\right)\cdot 3^k= b_0+
\sum\limits_{k=1}^\infty b_k 3^k,$$ shows that $x$ is a solution
of the equation \eqref{eq1.1}. \end{proof}

From the proof of Theorem \ref{t8} it is easy to see that if
$A_0=x^2_0+a_0 \equiv 0 \,(mod \, 3),$ then we have the following
congruences and appropriate equalities

a) $x^3_0 \equiv b_0 \,(mod \,3),$ i.e. $x^3_0 = b_0 +
M_1(x_0)\cdot 3 ;$

\

b) $x_0a_0+ M_1(x_0) \equiv b_1 \, (mod \,3),$ then $x_0a_0+
M_1(x_0) = b_1 + M_2(x_0)\cdot 3 ;$

\

c) $x_0a_1 + M_2(x_0)\equiv b_2 \,(mod \,3),$ then
\begin{equation}\label{eq3.4}x_0a_1 + M_2(x_0) =  b_2 + M_3(x_0)\cdot 3;\end{equation}

d) $\frac {A_0} 3 x_1 + x_1a_1 +
x_0a_2+x_0x_1^2+x_1^3+M_3(x_0)\equiv b_3 \,(mod \,3),$ it follows
that
$$\frac {A_0} 3 x_1 + x_1a_1 +
x_0a_2+x_0x_1^2+x_1^3+M_3(x_0) =  b_3 + M_4(x_0,x_1)\cdot 3.$$

Since $A_1 = \frac {A_0} 3 +a_1+2x_0x_1,$ then the congruence d)
can be written in the form
$(A_1-x_0x_1)x_1+x_0a_2+x_1^3+M_3(x_0)\equiv b_3 \,(mod \,3),$ and
 so we have
$$(A_1-x_0x_1)x_1+x_0a_2+x_1^3+M_3(x_0) =  b_3 + M_4(x_0,x_1)\cdot 3.$$

If for any natural number $k$ we have $A_k\equiv 0\, (mod\, 3),$
then we could establish the criteria of solvability for the
equation \eqref{eq1.1}. However, if there exists $k$, such that
$A_k\not\equiv 0 \, (mod \,3),$ then the criteria of solvability
can be found, and therefore, we need the following

\begin{lem}\label{l2} Let $\gamma(a)=1,$ $\gamma(b)=0$ and $x\in
\mathbb{Z}^{*}_3$ to be so that $A_{k-j}\equiv 0 \, (mod \,3), 1
\leq j \leq k,$ $A_k\not\equiv 0 \, (mod \,3)$ for some fixed $k$.
If $x$ be a solution of the equation \eqref{eq1.1}, then it is
true the following system of the congruences
\begin{equation}\label{eq3.5}
\begin{array}{c}
x^3_0 \equiv b_0 \, (mod \, 3),\\
x_0a_0+M_1(x_0) \equiv b_1 \, (mod \, 3),\\
x_{j-1}a_j+x_{j-2}a_{j+1}+\dots+x_0a_{2j-1}
+S_{2j}^{j}+M_{2j}(x_0,x_1,\dots,x_{j-1})\equiv b_{2j} (mod
\,3),\\
(A_j-x_0x_j)x_j+x_{j-1}a_{j+1}+x_{j-2}a_{j+2}+\dots+x_0a_{2j}
+S_{2j+1}^{j}+\\ +M_{2j+1}(x_0,x_1,\dots,x_{j-1})\equiv b_{2j+1}
(mod \,3),\\
A_kx_{k+i}+x_{k+i-1}a_{k+1}+x_{k+i-2}a_{k+2}+\dots+x_0a_{2k+i}
+S_{2k+1+i}^{k+i}+ \\ +M_{2k+1+i}(x_0,x_1,\dots,x_{k+i-1})\equiv
b_{2k+1+i} (mod \,3),\end{array}\end{equation} where $1 \leq j
\leq k$ and integers $M_{k}(x_0,\dots,x_{k-2})$ are defined from
the equalities
\begin{equation}\label{eq3.6}
\begin{array}{c}3\cdot M_1(x_0)=x^3_0 -b_0,\\ 3\cdot M_2(x_0) = x_0a_0+M_1(x_0) - b_1,\\
3\cdot M_{2j+1}(x_0,\dots,x_{j-1}) = x_{j-1}a_j+\dots+x_0a_{2j-1}
+S_{2j}^{j}+M_{2j}(x_0,\dots,x_{j-1}) -  b_{2j}, \\ 3\cdot
M_{2j+2}(x_0,x_1,\dots,x_j)
=(A_j-x_0x_j)x_j+x_{j-1}a_{j+1}+x_{j-2}a_{j+2}+ \dots+x_0a_{2j}+\\
+S_{2j+1}^{j}+M_{2j+1}(x_0,x_1,\dots,x_{j-1})- b_{2j+1},\\ 3\cdot
M_{2k+2+i}(x_0,x_1,\dots,x_{k+i}) = A_kx_{k+i}+x_{k+i-1}a_{k+1}
+x_{k+i-2}a_{k+2}+\dots+x_0a_{2k+i}+\\
+S_{2k+1+i}^{k+i}+M_{2k+1+i}(x_0,x_1,\dots,x_{k+i-1})
-b_{2k+1+i}.\end{array}\end{equation}
\end{lem}

\begin{proof} We will prove Theorem by induction. Let $k=1,$ i.e.
$$A_0=x^2_0+a_0 \equiv 0 \, (mod \, 3), \quad A_1=\displaystyle\frac{A_0}{3}+a_1+2x_0x_1\not\equiv 0 \, (mod \,
3),$$ then the system of the congruences \eqref{eq3.4} are true.
Note that $S_2^1 = 0,$ $S_3^1 = x_1^3.$

From \eqref{eq3.1} it is easy to get
$$\frac {A_0} 3 x_{t-2}+x_{t-2}a_1+x_{t-3}a_2+\dots+x_0a_{t-1}+ S_t^{t-2}
+2x_0x_1x_{t-2}+$$$$+M_t(x_0,\dots,x_{t-3}) \equiv b_t
 \, (mod \,p), \quad t\geq 4.$$ Therefore,
\begin{equation}\label{eq3.7}A_1 x_{t-2}+x_{t-3}a_2+\dots+x_0a_{t-1}+ S_t^{t-2}
+M_t(x_0,\dots,x_{t-3}) \equiv b_t
 \, (mod \,3), \quad t\geq 4,\end{equation} where
$$3\cdot M_{t+1}(x_0,\dots,x_{t-2}) = A_1 x_{t-2}+x_{t-3}a_2+\dots+x_0a_{t-1}+ S_t^{t-2}
+M_t(x_0,\dots,x_{t-3}) - b_t, \quad t\geq 4. $$

Obviously, the statement of Lemma is true for $k=1,$ i.e. for $i=t-3.$

Let $k=2,$ i.e. $A_0\equiv 0 \, (mod \,3),$ $A_1\equiv 0 \, (mod
\,3)$ и $A_2=\displaystyle\frac{A_1}{3}+a_2+x^2_1+2x_0x_2
\not\equiv 0 \, (mod \, 3),$ then from the equalities
\eqref{eq3.7} it follows that the following congruences are be
added to the system \eqref{eq3.4}:

e) $x_1a_2+x_0a_3+S_4^2 +M_4(x_0,x_1)\equiv b_4 \,(mod \,3),$ it
follows $$3 \cdot M_5(x_0,x_1)=x_1a_2+x_0a_3+S_4^2 +M_4(x_0,x_1)-
b_4;$$

f) $\frac {A_1} 3x_2+ x_2a_2+x_1a_3+x_0a_4+S_5^3
+M_5(x_0,x_1)\equiv b_5 \,(mod \,3),$ it follows $$3 \cdot
M_6(x_0,x_1,x_2)=\frac {A_1} 3x_2+ x_2a_2+x_1a_3+x_0a_4+S_5^3
+M_5(x_0,x_1)- b_5;$$
$$h) \frac {A_1} 3x_{t-2}+
x_{t-2}a_2+x_{t-3}a_3+\dots+x_0a_t+S_{t+1}^{t-1} +M_{t+1}(x_0,x_1,
\dots, x_{t-3})\equiv b_{t+1} \,(mod \,3),$$

where $t \geq 5$ and $M_{t+2}(x_0,x_1, \dots, x_{t-2})$ are
defined by equalities
$$3 \cdot M_{t+2}(x_0,x_1, \dots, x_{t-2}) = \frac {A_1} 3x_{t-2}+
x_{t-2}a_2+x_{t-3}a_3+\dots+x_0a_t+$$$$+S_{t+1}^{t-1}
+M_{t+1}(x_0,x_1, \dots, x_{t-3})- b_{t+1}.$$

Since $S_4^2=0,$ $S_2^2=0,$ $S_5^3=x_0x_2^2+x_1^2x_2,$
$S_{t+1}^{t-1}=S_{t+1}^{t-2}+x_1^2x_{t-2}+2x_0x_2x_{t-2},$ we
denote by $i = t-4$ and have

e) $x_1a_2+x_0a_3 +M_4(x_0,x_1)\equiv b_4 \,(mod \,3),$

\

f) $(A_2 - x_0x_2)x_2+ x_1a_3+x_0a_4+ M_5(x_0,x_1)\equiv b_5
\,(mod \,3),$

\

h)
$A_2x_{i+2}+x_{i+1}a_3+x_{i}a_4+\dots+x_0a_{i+4}+S_{i+5}^{i+2}+M_{i+5}(x_0,x_1,
\dots, x_{i+1})\equiv b_{i+5} \,(mod \,3),$ where $$3 \cdot
M_5(x_0,x_1)=x_1a_2+x_0a_3+M_4(x_0,x_1)- b_4,$$
$$3 \cdot M_6(x_0,x_1,x_2)=(A_2 - x_0x_2)x_2+x_1a_3+x_0a_4+M_5(x_0,x_1)-b_5.$$
$$3 \cdot M_{6+i}(x_0, \dots, x_{i+2}) = A_2x_{i+2}+
x_{i+1}a_3+\dots+x_0a_{i+4}+S_{i+5}^{i+2}+M_{i+5}(x_0,\dots,
x_{i+1})- b_{i+5}.$$

So we showed that the statement of Lemma is true for $k=2.$

Let the system of congruences \eqref{eq3.5}-\eqref{eq3.6} is true
for $k.$ Since $A_k\equiv 0 (mod \,3),$ then from the congruences
$$A_kx_{k+i}+x_{k+i-1}a_{k+1}+x_{k+i-2}a_{k+2}+\dots+x_0a_{2k+i} +S_{2k+1+i}^{k+i}+
$$$$+M_{2k+1+i}(x_0,x_1,\dots,x_{k+i-1})\equiv b_{2k+1+i} (mod \,3),\, i\geq 1$$
we derive
$$x_ka_{k+1}+x_{k-1}a_{k+2}+\dots+x_0a_{2k+1} +S_{2k+2}^{k+1}
+M_{2k+2}(x_0,\dots,x_k)\equiv b_{2k+2} (mod \,3),$$
$$\frac {A_k} 3 x_{k+1}+x_{k+1}a_{k+1}+x_{k}a_{k+2}+\dots+x_0a_{2k+2} +S_{2k+3}^{k+2}+
M_{2k+3}(x_0,\dots,x_{k})\equiv b_{2k+3} (mod \,3),$$
$$\frac {A_k} 3 x_{k+1+i}+x_{k+1+i}a_{k+1}+x_{k+i}a_{k+2}+\dots+x_0a_{2k+i+2} +S_{2k+i+3}^{k+i+2}+
$$$$+M_{2k+i+3}(x_0,x_1,\dots,x_{k+i})\equiv b_{2k+i+3} (mod \,3), \quad i \geq 1.$$

It is easy to check that $$S_{2k+3}^{k+2} = S_{2k+3}^{k+1} +
R_{k+1}x_{k+1} -x_0x_{k+1}^2,\qquad S_{2k+i+3}^{k+i+2} =
S_{2k+i+3}^{k+i+1} + R_{k+1}x_{k+1+i}, \quad i \geq 1.$$

By these correlations we deduce
$$\frac {A_k} 3 x_{k+1}+x_{k+1}a_{k+1}+x_{k}a_{k+2}+\dots+x_0a_{2k+2} +S_{2k+3}^{k+2}+
M_{2k+3}(x_0,x_1,\dots,x_{k}) = $$
$$=(A_{k+1}  -x_0 x_{k+1})x_{k+1}+x_{k}a_{k+2}+\dots+x_0a_{2k+2} +
S_{2k+3}^{k+1} + M_{2k+3}(x_0,x_1,\dots,x_{k}).$$

For $i \geq 1$ we get
$$\frac {A_k} 3 x_{k+1+i}+x_{k+1+i}a_{k+1}+\dots+x_0a_{2k+i+2} +S_{2k+i+3}^{k+i+2}+
M_{2k+i+3}(x_0,x_1,\dots,x_{k+i}) = $$
$$=A_{k+1}x_{k+1+i}+x_{k+i}a_{k+2}+\dots+x_0a_{2k+i+2}
+S_{2k+i+3}^{k+i+1} + M_{2k+i+3}(x_0,x_1,\dots,x_{k+i}).$$

Consequently, we have
$$x_ka_{k+1}+x_{k-1}a_{k+2}+\dots+x_0a_{2k+1} +S_{2(k+1)}^{k+1}+
M_{2(k+1)}(x_0,\dots,x_k)\equiv b_{2(k+1)} (mod \,3),$$
$$(A_{k+1}  -x_0 x_{k+1})x_{k+1}+x_{k}a_{k+2}+\dots+x_0a_{2k+2} +
S_{2k+3}^{k+1}+M_{2k+3}(x_0,\dots,x_{k})\equiv b_{2k+3} (mod
\,3),$$
$$A_{k+1}x_{k+1+i}+x_{k+i}a_{k+2}+\dots+x_0a_{2k+i+2}
+S_{2k+i+3}^{k+i+1}+ $$$$+ M_{2k+i+3}(x_0,x_1,\dots,x_{k+i})\equiv
b_{2k+i+3} (mod \,3), \quad i\geq 1 .$$

So we established that the system of congruences
\eqref{eq3.5}-\eqref{eq3.6} is true for $k+1.$ \end{proof}

Using the Lemma \ref{l2} we obtain the following Theorems.

\begin{thm}\label{t9} Let $\gamma(a)=1,$ $\gamma(b)=0$ and $x\in
\mathbb{Z}^{*}_3$ to be such that $A_{k-j}\equiv 0 \, (mod \,3), 1
\leq j \leq k,$ $A_k\not\equiv 0 \, (mod \,3)$ for some fixed $k\,
(k \geq 1)$. \emph{Then} $x$ to be a solution of the equation
\eqref{eq1.1} in $\mathbb{Z}^*_3$ if and only if the system of the
congruences
$$x^3_0 \equiv b_0 \, (mod \, 3),$$
$$x_0a_0+M_1(x_0) \equiv b_1 \, (mod \, 3),$$
$$x_{j-1}a_j+x_{j-2}a_{j+1}+\dots+x_0a_{2j-1}
+S_{2j}^{j}+M_{2j}(x_0,x_1,\dots,x_{j-1})\equiv b_{2j} (mod
\,3),$$
$$(A_j-x_0x_j)x_j+x_{j-1}a_{j+1}+x_{j-2}a_{j+2}+\dots+x_0a_{2j} +S_{2j+1}^{j}+$$
$$+M_{2j+1}(x_0,x_1,\dots,x_{j-1})\equiv b_{2j+1} (mod \,3),$$
has a solution, where $1 \leq j \leq k$ and integers
$M_{k}(x_0,x_1,\dots,x_{k-1})$ are defined from the equalities
$$3\cdot M_1(x_0)=x^3_0 -b_0,$$
$$3\cdot M_2(x_0) = x_0a_0+M_1(x_0) - b_1,$$
$$3\cdot M_{2j+1}(x_0,\dots,x_{j-1}) = x_{j-1}a_j+x_{j-2}a_{j+1}+\dots+x_0a_{2j-1}
+S_{2j}^{j}+M_{2j}(x_0,\dots,x_{j-1}) -  b_{2j},$$
$$3\cdot M_{2j+2}(x_0,x_1,\dots,x_j) =(A_j-x_0x_j)x_j+x_{j-1}a_{j+1}+x_{j-2}a_{j+2}+
\dots+x_0a_{2j}+$$
$$+S_{2j+1}^{j}+M_{2j+1}(x_0,x_1,\dots,x_{j-1})- b_{2j+1}.$$
\end{thm}

\begin{thm}\label{t10} Let $\gamma(a)=1,$ $\gamma(b)=0$ and $x\in
\mathbb{Z}^{*}_3$ to be so that $A_k\equiv 0 \, (mod \,3)$ for all
$k\in \mathbb{N}.$ Then $x$ to be a solution of the equation
\eqref{eq1.1} in $\mathbb{Z}^*_3$ if and only if the system of the
congruences
$$x^3_0 \equiv b_0 \, (mod \, 3),$$
$$x_0a_0+M_1(x_0) \equiv b_1 \, (mod \, 3),$$
$$x_{j-1}a_j+x_{j-2}a_{j+1}+\dots+x_0a_{2j-1}
+S_{2j}^{j}+M_{2j}(x_0,x_1,\dots,x_{j-1})\equiv b_{2j} (mod
\,3),$$
$$(A_j-x_0x_j)x_j+x_{j-1}a_{j+1}+x_{j-2}a_{j+2}+\dots+x_0a_{2j} +S_{2j+1}^{j}+$$
$$+M_{2j+1}(x_0,x_1,\dots,x_{j-1})\equiv b_{2j+1} (mod \,3),$$
has a solution, where $j\geq 1$ and integers
$M_{k}(x_0,x_1,\dots,x_{k-1})$ are defined from the equalities
$$3\cdot M_1(x_0)=x^3_0 -b_0,$$
$$3\cdot M_2(x_0) = x_0a_0+M_1(x_0) - b_1,$$
$$3\cdot M_{2j+1}(x_0,\dots,x_{j-1}) = x_{j-1}a_j+x_{j-2}a_{j+1}+\dots+x_0a_{2j-1}
+S_{2j}^{j}+M_{2j}(x_0,\dots,x_{j-1}) -  b_{2j},$$
$$3\cdot M_{2j+2}(x_0,x_1,\dots,x_j) =(A_j-x_0x_j)x_j+x_{j-1}a_{j+1}+x_{j-2}a_{j+2}+
\dots+x_0a_{2j}+$$
$$+S_{2j+1}^{j}+M_{2j+1}(x_0,x_1,\dots,x_{j-1})- b_{2j+1}.$$
\end{thm}

\section{Acknowledgement}

The first author was supported by grant UniKL/IRPS/str11061, Universiti Kuala Lumpur.

\

{\sc Rikhsiboev I.M.} Universiti Kuala Lumpur - (UniKL) Malaysian
Institute of Industrial Technology - (MITEC), e-mail:
\emph{ikromr@gmail.com}

{\sc Khudoyberdiyev A.Kh., Kurbanbaev T.K., Masutova K.K.}
Institute of Mathematics, 29, Do'rmon yo'li srt., 100125, Tashkent
(Uzbekistan) e-mails: \emph{khabror@mail.ru --- tuuelbay@mail.ru
--- kamilyam81@mail.ru}

\end{document}